\documentclass[a4paper,notitlepage, 8pt,reqno]{article}
  \usepackage{graphicx}
  \usepackage{mathtext}
 \usepackage[cp1251]{inputenc}
  \usepackage[T2A]{fontenc}
 \usepackage[russian,english]{babel}
 \usepackage[all,arc,poly,2cell,curve,arrow,tips]{xy}
  \usepackage{enumerate, eucal, amsthm,amsmath, amssymb}

\usepackage{tikz}
\usetikzlibrary{positioning}

  \allowdisplaybreaks[4]

 \theoremstyle{definition}
 \newtheorem{defn}{Definition}

 \theoremstyle{plain}
 \newtheorem{thm}{Theorem}
 \newtheorem*{thm*}{Theorem}
 
  \newtheorem*{prop*}{Proposition}
 
  \newtheorem*{cor*}{Следствие}
 
  \newtheorem*{lem*}{Лемма}

 \theoremstyle{remark}

 \newtheorem*{remark*}{Замечание}

  \newcounter{ab}
\title{The Stokes phenomenon for an irregular Gelfand-Kapranov-Zelevinsky system associated with the rank one lattice}

 \author{D. V. Artamonov\footnote{artamonov.dmitri@gmail.com}}

\begin{document}
\maketitle

\begin{abstract}
An explicit description of a multidimensional Stokes phenomenon for a Gelfand-Kapranov-Zelevinsky  system associated with a lattice of rank one is given.
\end{abstract}

\section{Introduction}

A keystone in the analytic theory of differential equations is  a construction called the Riemann-Hilbert correspondence. It is correspondence between systems of linear differential equations and monodromy data, which characterize a behaviour of solutions of the system near singular points of the system. The construction of the monodromy data differs in the case of a regular singular point and in the case of an irregular  singular point. In the case of the irregular point the monodromy is not enough to describe uniquely the behaviour os solutions, it is necessary to add Stokes matrices, that describe the Stokes phenomenon in the irregular singular point (see Section \ref{Stox}).

In the multidimensional analog of this theory one consider  an integrable in the Frobenious sense connection in holomorphic vector bundles. Locally they are presented as partial differential equations of special type - pfaffian systems. Also one considers even more general object - holonomic $\mathcal{D}$-modules.  The construction of generalized monodromy data  on this language in dimension $1$ is given in Section \ref{dmod1}.

As in dimension $1$ one considers  the regular and irregular cases. The construction of the Riemann-Hilbert correspondence in the regular case is straightforward, but the construction in the irregular case is much more difficult.
Some important elements of this construction go back to \cite{D}, intermediate results were obtained in \cite{Mal},
\cite{Mal91}, \cite{BV}, \cite{S0}, the final construction was given in \cite{M} (see also \cite{S}, where the systematic explanation is given). This construction is shortly given in Section \ref{stoxn}.

Even in dimension one not so many examples of explicit description of the Stokes phenomenon in the irregular case are known. Despite some artificial examples (see \cite{Z}) one can list the Jordan-Pochgammer equation
 \cite{In}, a parabolic cylinder equation,  an equation for the Eiri functions \cite{I}, an equation for the generalized hypergeometric functions, associated with a series $F_{p,q}$
 \cite{Du}.  These  hypergeometric functions, equation for them, the Stokes phenomenon for this equation is described in Section \ref{fpqs}.

 Mention that in dimension greater that one there appears a new phenomenon. To construct a Riemann-Hilbert correspondence one has to consider only system that have good formal decomposition. Also only for these system there exists a description of generalized monodromy data analogous to  one-dimensional case.

At the same time  there are no  examples of explicit description of Stokes phenomenon in the multidimensional case. In the present paper we give such a description for an irregular Gelfand-Kapranov-Zelevinsky (GKZ for short) system
(see \cite{Gel}), associated with a lattice of rank $1$ (see Section
\ref{sgkz}). The solutions of such system can be expressed through the generalized hypergeometric functions of one argument associated with the  series $F_{p,q}$.

Using this relation we manage to describe completely the local Riemann-Hilbert correspondence for the GKZ system.  This description is the main result of the paper, see Section
\ref{stoxgkz}.  In Section \ref{rhf}  a local Riemann-Hilbert functor is defined, in Section \ref{stoxgkz} we compute the value of this functor for the considered GKZ system.

Note that there are not so many papers in which irregular GKZ system is considered.  Mention only the papers \cite{f1}, \cite{f2}, where an irregularity sheaf of some GKZ system(a sheaf, which as multidimensional generalization of a Katz index in the case of dimension $1$).

\section{The Stokes phenomenon in dimension $1$} \label{Stox}

Let us  explain what is the Stokes phenomenon in the case of ordinary differential equations, let us define new objects are added to the monodromy (see \cite{I}, \cite{Sib})

\subsection{An approach based on formal normalization}

Let us be given a system

\begin{equation}
\label{osneq} \frac{dy}{dx}=A(x)y,
\end{equation}

where $y=(y_1,...,y_n)^t$ is a column-vector and $A(x)$ is an
$n\times n$-matrix.

\begin{defn} A point $x=0$ is an irregular singular point of the system \eqref{osneq}, if $0$ is a singularity of the matrix $A(x)$, and there exists a solution of   \eqref{osneq}, which has an exponential growth in $0$.
\end{defn}

A necessary condition for irregularity of the singularity $x=0$ is the existence of a pole of  $A(x)$ in zero of order greater then $1$. Let us write

\begin{equation*}
A(x)=(\sum_{k=-r-1}^{\infty}A_{k+1}x^{k})y
\end{equation*}

However the order of the pole $r+1$ can be changed under gauge transformations (i.e. linear transformations $y(z)\mapsto
B(z)y(z)$).  Suggest that applying a gauge transformation we made $r$ as minimal as possible.
Let the matrix $A_{-r}$ have different eigenvalues
$\alpha_1,...,\alpha_n$.  Then one says that $z=0$ is  {\it an unramified irregular singular point}.

Let zero be an  unramified irregular singular point. Let $P$ be a transformation that generalizes $A_{-r}$:
$$P^{-1}A_{-r}P=\Lambda_{-r}:=diag(\alpha_1,...,\alpha_n).$$

There exist the following theorem called the theorem  about {\it  formal normalization}.

\begin{thm}
The equation \eqref{osneq} in a neighborhood of an unramified
irregular singular point has a unique formal fundamental solution of
type
\begin{equation}
\label{fo}
Y=P(\sum_{k=0}^{\infty}Y_kx^k)exp(\frac{\Lambda_{-r}}{-r}x^{-r}+...+\frac{\Lambda_{-1}}{-1}x^{-1}+\Lambda_0lnx),
\end{equation}

where matrices $\Lambda_i$ are diagonal, $Y_0=I$ and $\sum_{k=0}^{\infty}Y_kx^k$  is a formal power series.
\end{thm}
In other words, the equation \eqref{osneq} has a formal fundamental solution of type
\eqref{fo}.

Also there exists a theorem about  {\it analytic sectorial normalization }. To give it's formulation let us define {\it  Stokes lines}
$l_{i,j}$, $i,j=1,...,n$ by formula

\begin{equation*}
l_{i,j}=\{x: arg \frac{\alpha_i-\alpha_j}{r}+\frac{2\pi
k}{r}\},\,\,\, k=0,...,r-1.
\end{equation*}

Also Stokes lines can be defined as follows. Put
$\Lambda(x)=\frac{\Lambda{-r}}{-r}x^{-r}+...+\frac{\Lambda_{-1}}{-1}x^{-1}=diag(\Lambda_1(x),...,\Lambda_n(x))$.
To the angle $\varphi$ there corresponds a ray that goes from $0$ in the direction $0$
$\varphi$. We say that $\Lambda_i(x)<_{\theta}\Lambda_j(x)$, if
$e^{\Lambda_i(x)-\Lambda_j(x)}$ has a power-like growth if one approaches to $0$ in the direction $\varphi$. Then Stokes directions for a pair
$i,j$  are such directions that one has neither $\Lambda_i(x)<_{\theta}\Lambda_j(x)$, nor
$\Lambda_j(x)<_{\theta}\Lambda_i(x)$.

Let us define{\it a Stokes sector} $S$   as a maximal sector $S=\{x:\,\theta_1<arg
x<\theta_2\}$, such that it does not contain two Stokes lines for any pair $i,j$.

\begin{thm}
Inside a Stokes sector  $S$  there exists an analytic fundamental solution  $Y_S$ that has the formal solution \eqref{fo} as an asymptotic expansion.
\end{thm}

However if one passes from one sector to another the analytic transformation  $P$ changes.  This surprising fact is called {\it the Stokes phenomenon}.

\begin{defn} Let the us be given the Stokes sectors $S_1$ and $S_2$  which have nonempty intersection. A {\it Stokes matrix} is defined by formula

\begin{equation}
\label{stot} St_{1,2}=Y_{S_1}Y_{S_2}^{-1}
\end{equation}

\end{defn}

\begin{defn}
In the case of an unramified irregular singular point we define it's {\it generalized monodromy data} as a collection of
a monodromy matrix, a collection of Stokes sectors and Stokes matrices.
\end{defn}

In the case of ramified singular point one must find a new variable $t$ instead of
  $x$ , such that they are related by formula

\begin{equation}
\label{zam}
 t^{d}=x,,\,\,\, d\in\mathbb{Z}_{+},
\end{equation}

and such that the point  $t=0$  becomes an unramified singular rpoint.

\subsection{An approach based on the Levelett-Hukuhara-Turrittin decomposition}

Another approach to the definition of Stokes matrices uses instead of formal normalization the
 {\it Levelett-Hukuhara-Turrittin decomposition}.

\begin{thm}
\label{to} Let $x=0$  be an uramified irregular singular point. There exists a formal transformation $y(x)\mapsto
P(x)y(x)$, such that the matrix  $A(x)$ is transformed to block-diagonal $diag(A_1(x),...,A_p(x))$, where the blocks are of type $A_{j}(x)=diag(A_{j,1}(x),...,A_{j,m_j}(x))$, where the block
$A_{j,k}(x)$ is a Jordan cell with the eigenvalue
$\lambda_j(x)$.
\end{thm}

In a Stokes sector $S$ there exists an analytic transformation $P_S(x)$ asymptotic to $P(x)$. Then the Stokes matrices are defined by formulas analogous to
\eqref{stot}, as a relation of two analytic transformations that do the Levelett-Hukuhara-Turittin decomposition in two Stokes sectors.

\section{The Stokes phenomenon on the language of $\mathcal{D}$-modules in dimension $1$} \label{dmod1}

Sometimes it is convenient to consider instead of linear differential equations holonomic    $\mathcal{D}$-modules (see \cite{S}).

\subsection{$\mathcal{D}$-modules}  To the system of differential equations  \eqref{osneq} there corresponds a connection in a trivial bundle
 that is defined by formula
\begin{equation*}
\nabla=d-A(x).
\end{equation*}
To this bundle with a connection there corresponds a  sheaf $M$ of $\mathcal{D}$-module formed by sections of the bundle,  the operator $\frac{d}{dx}$ acts as $\nabla$.

Let $\Omega^{.}$ be a de Rham complex, define the {\it de Rham functor}

\begin{equation*}
DR(M)=M\otimes \Omega^{.}
\end{equation*}

Then $$\mathcal{L}=\mathcal{H}^0(DR(M))$$ is a sheaf of solutions of the system
\eqref{osneq}, and
$H^0(DR(M))=\Gamma(H^0,\mathcal{H}^0(DR(M))$  a space of solutions of this system.

An arbitrary holonomic  $\mathcal{D}$-module on a space  $X$ being
restricted to a dense open subset is isomorphic to a
$\mathcal{D}$-module of sections of bundle with a connection (see
definition of the holonomicity in \cite{dmo}).

\subsection{A blowup $\widetilde{X}$ of a disk in zero and sheaves of holomorphic functions on it }

Let $X$ be  a small disk with a center in zero, denote as
$\widetilde{X}$ a real blowup of   $X$ in zero, also denote as
$\partial \widetilde{X}=S^1$ a boundary circle of this blowup. Let

\begin{equation*}
\omega: \widetilde{X}\rightarrow X
\end{equation*}

be a natural projection, hence $\partial
\widetilde{X}=\omega^{-1}(0)$. Put
$$X^*:=X\setminus\{0\},$$  note that from one hand  $X^*\subset
X$ and from the other hand  $X^*\subset \widetilde{X}$.

Define a sheaf $\mathcal{A}$ of holomorphic functions on
$\widetilde{X}$, and  sheaves  $\mathcal{A}^{mod}$  and
$\mathcal{A}^{rd}$ of holomorphic functions on $X^*$ that have moderate growth of random decay when one approaches to
к $\partial{\widetilde{X}}$. Define the sheaf
$e^{\varphi}\mathcal{A}^{mod}$ of holomorphic functions, that grow  not faster then $e^{\varphi}$ when one approaches to $\partial\widetilde{X}$.

\subsubsection{The sheaf of holomorphic functions}
Defined the operator $\bar{\partial}$ on $C_{\infty}(\widetilde{X})$.
Take polar coordinates
$r,\varphi$. In these coordinates

$$\bar{\partial}=\frac{1}{2}(r\partial_r+i\partial_{\varphi}).$$

Using this formula one can define the action of
$\bar{\partial}$  on  $C_{\infty}(\widetilde{X})$ and a sheaf  $\mathcal{A}$ of holomorphic functions on $\widetilde{X}$:

\begin{equation}
\mathcal{A}:=Ker\overline{\partial}.
\end{equation}

\subsubsection{Sheaves of functions of moderate growth and of random decay on $\widetilde{X}$ } Let us define the sheaves $\mathcal{A}^{mod}$  and
$\mathcal{A}^{rd}$ of functions of moderate growth and of random decay on $\widetilde{X}$.

The first sheaf is defined as follows. For every open subset
$\widetilde{U}\subset\widetilde{X}$ a section $f$ of the sheaf
$\mathcal{A}^{mod}$ is a holomorphic function on
$U^*=\widetilde{U}\cap X^*$, such that for every compact subset
$K\subset \widetilde{U}$, such that in it's neighborhood $D$ is
defined by function $g_K\in\mathcal{O}_X(K)$, for some $C> 0$,
$N\geq 0$ (that depend on
 $K$) one has

$$|f|\leq C_K|g_K|^{-N_K}.$$

The sheaf $\mathcal{A}^{rd}$ is defined analogously but one  claims that for every and  $N$ and every $K$ there exists  $C>0$, such that

$$|f|\leq C_{K,N} |g_K|^N.$$

\subsubsection{The sheaf $e^{\varphi}\mathcal{A}^{mod}$} Define a sheaf
$e^{\varphi}\mathcal{A}^{mod}$  of functions of exponential growth not faster than $\varphi$, where $\varphi\in \Gamma(U,\mathcal{O}(X^*))$.

It's restriction onto the direction $\theta$ is defined by equality

$$(e^{\varphi}\mathcal{A}^{mod})_{\theta}=e^{\varphi_{\theta}}\mathcal{A}_{\theta}^{mod},$$
where $\varphi_{\theta}$ and $\mathcal{A}_{\theta}^{mod}$ is a restriction of the section and of the sheaf onto the ray the starts from zero in the direction $\theta$. Note that this sheaf depends only on the class of $\varphi$ in $\mathcal{O}_{X^*}\setminus \mathcal{O}_{X}$.

\subsection{$\mathcal{J}$-filtration in dimension $1$ }

The information carried by the generalized monodromy data and the change of variables \eqref{zam},
that the removes the ramification is carried by a filtration on the sheaf $\mathcal{L}$.
Let us define it.

Let $\mathcal{P}$ be a ring of singular parts of Loran expansions in zero (i.e. the ring of polynomials in $\frac{1}{x}$ without  constant).  Denote as $\mathcal{J}_d$ a locally constant sheaf on  $S^1$ with a stalk $\mathcal{P}$ a the monodromy $\varphi(x)\mapsto
\varphi(e^{2\pi i/d}x)$ corresponding to a bypass around the circle.  Thus a section $\varphi$ can be identified with a polynomial in $\frac{1}{\sqrt[d]{x}}$ without constant. This identification defines an embedding $\mathcal{J}_d\subset \mathcal{J}_{d'}$,  in the case when $d$ divides  $d'$.

Put

\begin{equation}
\mathcal{J}:=\bigcup \mathcal{J}_d
\end{equation}

Define an order on the stalks of this sheaf. Let $\theta\in S^1$,
we say that $\varphi\leq_{\theta} \psi$, if $e^{\varphi-\psi}$ has a power-like growth in the direction  $\theta$.

We define a $\mathcal{J}$-filtration of the sheaf  $\mathcal{L}$ over $S^1$ a family of subsheaves $\mathcal{L}_{\leq \varphi}$ in
$\mathcal{L}$ indexed by $\varphi\in\mathcal{P}$, such that for a direction $\theta$ in the case  $\varphi\leq_{\theta} \psi$ one has an inclusion of restrictions
$\mathcal{L}_{\leq \varphi,\theta}\subset \mathcal{L}_{\leq
\psi,\theta}$

\subsection{A filtration on the sheaf of solution}

In \cite{S} the following filtration on  $DR(M)$ is defined:

\begin{equation}
DR_{\leq \varphi}(M):=DR(e^{\varphi}\mathcal{A}^{mod}\otimes M).
\end{equation}

In other words $DR_{\leq \varphi}(M)$ is a subsheaf in $DR(M)$, such that it's stalk over the direction $\theta$ is formed by solutions that grow in this direction not faster then $e^{\varphi}$ when one approaches to zero.

 This filtration induces a filtration on the sheaf of solutions $\mathcal{L}$:

 \begin{equation}
 \mathcal{L}_{\leq \varphi}=\mathcal{H}^0(DR_{\leq \varphi}(M)).
 \end{equation}

In other words $ \mathcal{L}_{\leq \varphi}$ is a subsheaf in $DR(M)$, such that it's stalk over the direction $\varphi$ is formed by solutions that grow in this direction not faster then $e^{\varphi}$ when one approaches to zero.

This filtration is graded with a finite set of exponential factors. The set of exponential factors is a set
$\Phi=\{c_1,...,c_q\}$ of elements of $\mathcal{P}$ such that

\begin{equation}\label{rj}
\mathcal{L}_{\leq \varphi}=\sum_{\psi\in\Phi}\beta_{\psi\leq
\varphi}gr_{\psi}\mathcal{L},
\end{equation}

where $gr_{\psi}\mathcal{L}$ is a locally constant sheaf and
$\beta_{\psi\leq \varphi}$ is an functor which restricts a sheaf
$\mathcal{B}$  to the subspace
$\{\theta: \varphi<_{\theta} \psi\}$  and expands it to $S^1\setminus \{\theta: \varphi<_{\theta}
\psi\} $ by zero.

\subsection{The Levelett-Turittin decomposition} The set of exponential factors can be found using an analog for $\mathcal{D}$-modules of the theorem \ref{to}.

Let $\widehat{\mathcal{O}}$ be a ring of formal power series in zero on $X$, let $E^{\varphi}$  be a $\mathcal{D}$-module of sections of a linear bundle with a connection $d+d\varphi$, where $\varphi$ is a section
$\mathcal{J}$.
\begin{thm}
There exists an isomorphism

\begin{equation}
M\otimes_{\mathcal{O}} \widehat{\mathcal{O}}\simeq
\bigoplus_{\phi\in\Phi}E^{\varphi}\otimes_{\mathcal{O}} R_{\varphi},
\end{equation}

where $R_{\varphi}$  is a  $\mathcal{D}$-module of sections of a bundle with a connection with only regular singularities and $\Phi$ is a finite set of section of
$\mathcal{J}$.
\end{thm}

\subsection{ A correspondence between filtrations and Stokes lines, sectors and matrices in dimension $1$ }

Let us show how one can reconstruct Stokes lines, sectors and matrices from the filtration.

Stokes lines are those directions $\theta$ for which some two exponential factors $\varphi_i,\varphi_j\in\Phi$  are incomparable. The sectors are defined as earlier.

It is important that on the Stoke sector the sheaves $\beta_{\varphi\leq \psi}gr_{\psi}\mathcal{L}$ are constant.

The Stokes matrices are defined as matrices connecting two filtration $\mathcal{L}_{\leq \varphi}$ in two sectors.  Let $I_1$,
$I_2$ be two Stokes sectors and let $\theta\in I_1\cap I_2$.

Consider the spaces $L_1=\Gamma(I_1,\mathcal{L})$ and
$L_2=\Gamma(I_2,\mathcal{L})$ and also  $L_1^{\psi}=\Gamma(I_1,gr_{\psi}\mathcal{L})$ и
$L_2^{\psi}=\Gamma(I_2,gr_{\psi}\mathcal{L})$. The filtration \eqref{rj} defines decompositions of these spaces into direct sums

$$L_1=\oplus_{\psi\in\Phi}L_1^{\varphi},\,\,\,L_2=\oplus_{\psi\in\Phi}L_2^{\psi}.$$

Also since the sheaves $\beta_{\varphi\leq \psi}gr_{\psi}\mathcal{L}$  are constant on  $I_1$  and $I_2$ one has isomorphisms

$$a_1: L_1\rightarrow \mathcal{L}_{\theta},,\,\,\, a_2: L_2\rightarrow \mathcal{L}_{\theta},$$

the the Stokes matrix is reconstructed as follows

$$S=a_2^{-1}a_1.$$

\section{ The Stokes phenomenon in dimension $> 1$} \label{stoxn}

To describe the Stokes phenomenon and the Riemann-Hilbert
correspondence in dimension greater then $1$ the language of
holonomic   $\mathcal{D}$-modules is used.  A $\mathcal{D}$-module
is called regular if it's restriction on every curve  is a
one-dimensional regular $\mathcal{D}$-module,  otherwise it is
called irregular.

\subsection{Pre-Stokes filtration } \label{prestox}

To describe the Stokes phenomenon one uses  filtration. Let $Y$ - be
a topological space and let $\mathcal{J}$ be a sheaf of ordered
abelian groups on it. Denote as $\mathcal{J}^{et}$  the etale space
of the sheaf   $\mathcal{J}$, introduce a notation
\begin{equation}
\mu:\mathcal{J}^{et}\rightarrow Y
\end{equation}
for a projection onto the base   $Y$. Consider the diagram

\begin{tikzpicture}
  \node (A) {$\mathcal{J}^{et}\times \mathcal{J}^{et} $};
  \node (B) [below=of A] {$\mathcal{J}^{et}$};
  \node (C) [right=of A] {$\mathcal{J}^{et}$};
  \node (D) [right=of B] {$Y$};
  \draw[-stealth] (A)-- node[left] {\small $p_2$} (B);
  \draw[-stealth] (B)-- node [below] {\small $\mu$} (D);
  \draw[-stealth] (A)-- node [above] {\small $p_1$} (C);
  \draw[-stealth] (C)-- node [right] {\small $\mu $} (D);
\end{tikzpicture}

Let $(\mathcal{J}^{et}\times \mathcal{J}^{et} )_{\leq}$ be a subsheaf formed by sections $(\varphi(x),\psi(x))$ over one point $x$ for which  $\varphi(x)\leq \psi(x)$.  Denote as  $\beta_{\leq}$ a functor from the category of sheave on  $\mathcal{J}^{et}\times
\mathcal{J}^{et} $ to itself that restricts a sheaf $\mathcal{F}$ onto $(\mathcal{J}^{et}\times \mathcal{J}^{et} )_{\leq}$,
and then continues it by zero on the rest part of  $\mathcal{J}^{et}\times
\mathcal{J}^{et} $.

\begin{defn}

 Define a  {\it pre-$\mathcal{J}$ filtration}  as a sheaf $\mathcal{F}_{\leq}$ on the space   $\mathcal{J}^{et}$ together with a morphism
  $\beta_{\leq}(p_1^{-1}\mathcal{F}_{\leq})\rightarrow p_2^{-1}\mathcal{F}_{\leq}$.
\end{defn}

Let us explain this definition.  One has an isomorphism of stalks

$$(p_1^{-1}\mathcal{L}_{\leq})_{\varphi(x),\psi(x)}=(\mathcal{L}_{\leq})_{\varphi(x)},\,\,\,\,
(p_2^{-1}\mathcal{L}_{\leq})_{(\varphi(x),\psi(x))}=(\mathcal{L}_{\leq})_{\psi(x)}.$$

Thus if for a pair ov sections one has
$\varphi(x) \leq \psi(x)$, then we have a homomorphism

$$(\mathcal{L}_{\leq})_{\varphi(x)}\rightarrow (\mathcal{L}_{\leq})_{\psi(x)}.$$

Denote a stalk of s pre-$\mathcal{J}$ filtration over a section $\varphi(x)$ as $\mathcal{L}_{\leq \varphi(x)}$.

\subsection{The space $\widetilde{X}$ }
\label{pvo}
Let us be given a $\mathcal{D}$-module on the space  $X$  with a singular divisor  $D=k_1D_1+...+k_pD_p$, which has only normal crossings.

To describe the Stokes phenomenon we need a space $\widetilde{X}$,  which is a real blow-up along the divisor
 $D$. Explicitly the space $\widetilde{X}$ can be described as follows:
\begin{equation}
\widetilde{X}=\bigoplus_{i=1}^p S^1L(D_i),
\end{equation}

where $L(D_i)$ - is a linear bundle, associated with the divisor
$D_i$ and $S^1L(D_i)$ are  spherical bundles associated with $L(D_i)$.


As before denote as $$\omega :\widetilde {X} \rightarrow
X$$ the natural projection. Also introduce a notation
$$\partial \widetilde{X}:=\omega^{-1}(D).$$

As before put $$X^*=X\setminus D.$$

Let

\begin{equation}
j:\partial X^* \rightarrow X,\,\,\,\,\widetilde{j}:X^* \rightarrow
\widetilde{X},
\end{equation}

also let

\begin{equation}
i: D\rightarrow X,\,\,\,\, \widetilde{i}: \partial
\widetilde{X}\rightarrow \widetilde{X}
\end{equation}
 be natural embeddings.

Analogously to the one-dimensional case one defines a sheaf
$\mathcal{A}^{mod}$  of moderate growth and a sheaf
$\mathcal{A}^{rd}$ of functions of random decay but now these sheaves are  restricted  to $\partial{\widetilde{X}}$.

\subsection{The sheaf $\mathcal{J}$}

Let us give a construction of a sheaf $\mathcal{J}$ on the space $Y=\partial
\widetilde{X}$. To describe the Stokes phenomenon we use a
 $\mathcal{J}$-filtration defined in Section \ref{prestox}.

Let  $D=k_1D_1+...+k_pD_p$ -- be a decomposition of  $D$ into
irreducible components, introduce  a notation
$d=(d_1,...,d_p)\in\mathbb{Z}_{+}^p$. Change if necessary $X$ to a
sufficiently small neighborhood of
 $D$, then one can suggest that there exists a  $d$-covering   $$\rho_d: X_d \rightarrow X$$ along
$D$ (that is along  $D_i$ this is a $d_i$-covering), and also its lifting $$\widetilde{\rho}_d: \widetilde{X}_d\rightarrow
\widetilde{X}.$$ Let $$\omega_d:\widetilde{\rho}_d \rightarrow
X_d$$  be a natural projection. Thus we have a diagram

\begin{tikzpicture}
  \node (A) {$\widetilde{X}_d$};
  \node (B) [below=of A] {$X_d$};
  \node (C) [right=of A] {$\widetilde{X}$};
  \node (D) [right=of B] {$X$};
  \draw[-stealth] (A)-- node[left] {\small  $\omega_d$ } (B);
  \draw[-stealth] (B)-- node [below] {\small $\rho_d$} (D);
  \draw[-stealth] (A)-- node [above] {\small $\widetilde{\rho}_d$} (C);
  \draw[-stealth] (C)-- node [right] {\small $\omega $} (D);
\end{tikzpicture}

The sheaf $\mathcal{J}$ is defined as follows.

Put

\begin{equation}
\label{jdt}
\widetilde{\mathcal{J}}_d=\widetilde{\rho}_{d,*}(\omega_{d,*}\mathcal{O}_{X_d}(*D))\cap
\widetilde{j}_{*}\mathcal{O}_{X^*}.
\end{equation}

This formula means the following,  $\mathcal{O}_{X_d}(*D)$ is a sheaf of meromorphic functions on  $X_d$ with singularity on  $D$,  the sheaf
$\omega_{d,*}\mathcal{O}_{X_d}(*D)$ is it's lifting onto
$\widetilde{X}_d$. When one applies $\widetilde{\rho}_{d,*}$ and takes an intersection with
 $\widetilde{j}_{*}\mathcal{O}_{X^*}$ one obtains a sheaf of  meromorphic functions in variables $\sqrt[d_1]{x_1},...,\sqrt[d_p]{x_p},x_{p+1},...,x_n$ on  $\widetilde{X}$ with singularity on $\partial \widetilde{X}$.


Now put
\begin{equation}
\mathcal{J}_d:=\widetilde{\mathcal{J}}_d/\widetilde{\mathcal{J}_d}\cap(j_*
\mathcal{O}(*D))^{lb},
\end{equation}

where $(j_* \mathcal{O}_{X^*})^{lb}$ is a sheaf of locally bounded  near
$D$ sections of the sheaf  $j_* \mathcal{O}_{X^*}$. Note that if a meromorphic function on $X$ is bounded along some direction starting on $D$,  then it is holomorphic on $X$.

Put

\begin{equation}
\widetilde{\mathcal{J}}:=\bigcup_d\widetilde{\mathcal{J}}_d,
\,\,\,\,\,\, \mathcal{J}:=\bigcup_d\mathcal{J}_d.
\end{equation}

To define an order on the stalks of these sheaves it is sufficient to define which section are less than zero. A subsheaf
 formed by these sections we denote as $\widetilde{\mathcal{J}}_{\leq 0}$,
put

\begin{equation}
\widetilde{\mathcal{J}}_{\leq 0}=\widetilde{\mathcal{J}}\cap
log\mathcal{A}^{mod D},
\end{equation}

where $log\mathcal{A}^{mod D} $ is a sheaf of logarithmic growth at $\partial\widetilde{X}$.  This order defines an order on
$\mathcal{J}$.

Note that $\mathcal{J}_1=\omega^{-1}(\mathcal{O}_X(*D)\setminus
\mathcal{O}_X)$.

\subsection{Good families of exponential factors}
A principle new phenomenon that occurs in the multidimensional case is an existence of a good formal decomposition. In the one-dimensional case such a decomposition exists automatically. The existence of such  decompositions is principle in the construction of Riemann-Hilbert correspondence. Also, only for systems  with such a decomposition one can construct multidimensional analogs of  Stokes lines, sectors and matrices. In this section we define the property of goodness and in the next Section we define a good formal decomposition.

 Let us be given $x\in X$ and as before let $D=k_1D_1+...+k_pD_p$ be a divisor with normal crossings. Let
 $x\in D_{i_1},...,D_{i_l}$, $x\notin
D_1,...,\widehat{D_{i_1}},...,\widehat{D_{i_l}},...,D_{i_{p}}$. Without less of generality one can suggest that  $x\in D_1,...,D_l$,
$x\notin D_{l+1},...,D_p$.

Consider first an unramified case. Choose coordinates in a neighbourhood of   $x$, such that this point has zero coordinated and $D_i=\{x_i=0\}$,
$i=1,...,l$.
\begin{defn}
 A germ $\eta\in \mathcal{O}_X(*D)\setminus \mathcal{O}_X$ is called {\it monomial}, if modulo $\mathcal{O}_X$ one has an equality
 $$\eta=x_1^{-m_1}...x_l^{-m_l}u(x_1,...,x_n),$$ where $m_i\in\mathbb{Z}^+$, $u\in\mathcal{O}_X$, $u(0,...,0,x_{p+1},...,x_n)
\neq 0$.
\end{defn}

 In the case when the germ $\eta$ is monomial one has

 \begin{align*}
 \eta\leq_{\theta} 0\Leftrightarrow \eta=0\text{ or } arg (u(0))-\sum_jm_j\theta\in(\frac{\pi}{2},\frac{3\pi}{2}) .
 \end{align*}

 \begin{defn}
 We say that the finite set $\Phi\subset \mathcal{O}_X(*D)\setminus\mathcal{O}_X$ is {\it good},
  if it consists of one element or for every pair  $\varphi,\psi \in \Phi$ their difference $\varphi-\psi$ is monomial.
 \end{defn}

Let us be given a good finite set $\Phi\subset
\mathcal{O}_X(*D)\setminus\mathcal{O}_X$,  for every pair
$\varphi\neq \psi$  one has
$\varphi-\psi=x_1^{-m_1}...x_l^{-m_l}u(x)$. Let us give a definition

\begin{defn}
{\it  A Stokes hypersurface } is a set
\begin{equation}\label{sth}
St(\varphi,\psi)=\{(\theta_1,...,\theta_l)\in (S^1)^l:\,\, \sum_j
m_j\theta_j-arg u(0)=\pm \frac{\pi}{2} mod \,\,2\pi\}.
\end{equation}
\end{defn}

In the ramified case the definitions are the same but one takes coordinates $t_i=x_i^{d_i}$ instead of $x_i$.

\subsection{The Levelett-Turittin decomposition}

As above we suggest that  $x\in D_1,...,D_l$, $x\notin
D_{l+1},...,D_p$.

Firstly suggest that we are in the unramified case.

Let
$I\subset \{1,...,p\}$ be a set of indices, introduce a notation

\begin{equation*}
D(I):=\sum_{i\in I}k_i D_i,\,\,\,D(I^c):=\sum_{i\notin I}k_i
D_i,\,\,\, D(I^0):=D(I)\setminus D(I^c).
\end{equation*}

Take $\varphi\in \mathcal{O}_{X}(*D)\setminus \mathcal{O}_{X}$.
Choose it's representative over a sufficiently small open subset
$U$.
 Denote as $\varphi_I$ a
restriction of $\varphi$ onto $U\setminus  D(I^c)$. Thus
$\varphi_I\in \mathcal{O}_{U}(*D)\setminus \mathcal{O}_{U}(*
D(I^c))$

Denote as  $\Phi_I$ an image of $\Phi$.

\begin{defn}
We say that  $M$  has {\it an unramified good formal decomposition}, if there exists a good subset $\Phi \subset
\mathcal{O}_{X}(*D)\setminus \mathcal{O}_{X}$, such that for every
$I\subset\{1,...,p\}$
 in a neighborhood $U$ of a points  $x$ one has a decomposition

\begin{align}
\label{hfr} \widehat{\mathcal{O}}_{U\cap \widehat{D(I^0)}
}\otimes_{\mathcal{O}_{U\setminus D(I^0)}} M\mid_{U\cap
D(I^0)}\simeq \sum_{\varphi_I\in\Phi_I} (E^{\varphi_I}\otimes_{\mathcal{O}_{U\setminus D(I^0)}}
R^{\varphi_I})\mid_{U\cap D(I^0)},
\end{align}
where $\widehat{\mathcal{O}}_{U\cap \widehat{D(I^0)}}$ is a completion of   $\widehat{\mathcal{O}}_{U}$ alogn $U\cap \widehat{D(I^0)}$
(i.e. formal power series in variables $x_{l+1},...,x_{p}$,
whose coefficients are holomorphic functions depending on $x_1,...,x_l$).
\end{defn}

Consider now a ramified case.

\begin{defn}
We say that $M$ has  {\it  a good formal decomposition }, if there exists a mapping $\rho_d: X_d\rightarrow X$, such that
$\rho_d^+M$ has an unramified good formal decomposition.
\end{defn}

Explicitly this definition means the following. Instead of claiming
an existence of the isomorphism \eqref{hfr}, we claim the existence
of $d=(d_1,...,d_p)$, such that

\begin{align}
 \widehat{\mathcal{O}}_{U_d\cap \widehat{D(I^0)}
}\otimes_{\mathcal{O}_{U_d\setminus D(I^0)}}
\rho^{+}_dM\mid_{U_d\cap D(I^0)}\simeq \sum_{\varphi_I\in\Phi_{d,I}}
(E^{\varphi_I}\otimes_{\mathcal{O}_{U\setminus D(I^0)}} R^{\varphi_I})\mid_{U_d\cap D(I^0)},
\end{align}

where $\Phi_{d,I}$ is obtained by restriction of a good subset
$\Phi_d\subset \mathcal{J}_d$, $U_d$ is an open subset in
$X^d$  that contains $\rho_d^{-1}(x)$. The sections
$\widehat{\mathcal{O}}_{U\cap \widehat{D(I^0)}}$ can be considered as formal power series in variables $x^{d_{l+1}}_{l+1},...,x^{d_p}_{p}$,
whose coefficients are holomorphic functions of
$x^{d_1}_1,...,x^{d_l}_l$.

 \subsection{The Riemann-Hilbert functor}
\label{rhf}

 Let us be given a $\mathcal{D}$-module of sections of a bundle with a connection with a singular divisor  $D$.
 Put

 \begin{equation}
\mathcal{L}:=\widetilde{i}^{-1}\widetilde{j}_*\mathcal{H}^0DR(M\mid_{X^*}).
 \end{equation}

That is the sheaf  $M$ is restricted to a subspace
$X^*=X\setminus D$, where it is nonsingular, then one takes a sheaf of it's flat sections $\mathcal{H}^0DR(M\mid_{X^*})$, then one applies the functor $\widetilde{i}^{-1}\widetilde{j}_*$ and one gets a sheaf on
$\partial \widetilde{X}$, which is described as follows.
Take an open subset  $U\subset \partial
\widetilde{X}$ such that $U=\partial \widetilde{X}\cap
\widetilde{U}$, where $ \widetilde{U}$ -
 is a sufficiently small subset in $\widetilde{X}$. Then a section $s$ of the sheaf $\mathcal{L}$ over $U$ is a section
 $\widetilde{s}$ of the sheaf  $\mathcal{H}^0DR(M\mid_{X^*})$ over $
\widetilde{U}$. The section $s$ and $w$ are equivalent if and only if $\widetilde{s}$ and $\widetilde{w}$ coincide in a neighbourhood of  $\partial \widetilde{X}$.

Also let us define (see notations in Section \ref{prestox})

\begin{equation}
M_{\mathcal{J}^{Et}}^{mod}:= \mu^{-1}(\widetilde{i}^{-1}(
\mathcal{A}^{mod}\otimes _{\omega^{-1}\mathcal{O}_X}
\omega^{-1}M)),\,\,\, DR_{\mathcal{J}^{et}}^{mod D}(M):=
M_{\mathcal{J}^{Et}}^{mod}\otimes \Omega^{.,mod}_{\mathcal{J}^{Et}}.
\end{equation}

Now let us define a $\mathcal{J}$-filtration

\begin{equation}
\mathcal{L}_{\leq}:=\mathcal{H}^0DR_{\mathcal{J}^{et}}^{mod D}(M),
\end{equation}


\begin{defn}
Define a Riemann-Hilber functor
$$RH(M):=(\mathcal{L},\mathcal{L}_{\leq}).$$
\end{defn}

Suggest that $M$ has a good formal decomposition. Let $x\in D_1,...,D_l$, $x\notin D_{l+1},...,D_p$.  Consider a circle  that intersects every Stokes hypersurface.
 Consider Stokes intervals which are  maximal intervals $\Theta$ that contain not greater that one intersection with every hypersurface $St(\varphi,\psi)$ for every
$\varphi\neq \psi\in \Phi$.


 Then for  $x$ in a neighborhood of this section and an arbitrary section  $\varphi(x)$ of  $\mathcal{J}$ over this point one gas

 \begin{equation}
 \label{razll}
\mathcal{L}_{\leq \varphi(x)}=\bigoplus_{\psi\in\Phi, \,\,\,
\psi(x)\leq \varphi(x)}gr_{\psi(x)}\mathcal{L}\mid_{\varphi(x)},
 \end{equation}
where $gr_{\psi(x)}\mathcal{L}$ are some local systems.


Hence a sheaf $\mathcal{L}_{\leq}$  is defined over $(D_1\cap...\cap D_l)\setminus (D_{l+1}\cup...\cup
 D_p)$ by a collection of hypersurfaces and Stokes matrices that connect two expansions

 \begin{equation}
\bigoplus_{\psi\in\Phi}\mathcal{L}\mid_{\varphi(x)}
 \end{equation}

on an intersection of neighborhoods of Stokes intervals.

Mention that a sheaf on $\mathcal{J}^{et}$ defined by the formula \eqref{razll} automatically is a pre-Stokes filtration.

\subsection{A relation to isomonodromic deformation}

Integrable pfaffian systems in several variables are closely
related to isomonodromic deformations of equations in one variable
(for example  of the equation obtained by fixation of all variables
except one).  This pfaffian system is called an  {\it integrable
deformation }  of an ordinary equation. In the case when this
multidimensional system is regular, the integrability is equivalent
to isomonodromy.

In the case considered in the present paper an ordinary equation has
irregular singularities. There exist several approaches to the
definition of what is an isomonodromic deformation of an equation
with an irregular singularity. In \cite{I} for example an
isomonodromic deformation is a pfaffian system such that when one
changes the positions of singularities

\begin{enumerate}
\item a Jordan type of the normal form of the coefficient at the singularity in the expansion of the matrix of the connection into a Loran series is constant
\item the Stokes lines and sectors at new positions of singularities are obtained by parallel transport form the Stokes lines and sectors of the initial position of singularity.
\end{enumerate}

In \cite{An} it is shown, what happens if one omits the second claim. Then at almost all positions of singularities an integrable deformation still is isomonodromic (i.e. the second claim takes place). The existence of a good formal decomposition provides that the second claim holds in all points.

\section{Generalized hypergeometric series $F_{p,q}$}
\label{fpqs}

A generalized hypergeometric series in one variable is defined by formula (see \cite{be}, \cite{lu}, \cite{jfl})

\begin{equation}
F_{p,q}(a_1,...,a_p;b_1,...,b_q;z)=\sum_{n}\frac{(a_1)_n...(a_p)_n}{(b_1)_n...(b_q)_nn!}x^n,
\end{equation}
where $(a)_n=a(a+1)...(a+n-1)$ for $n>0$ and $(a)_0=1$.

This series does not converge for all $p,q$. More precise this series converges for all $x$ when $p\leq q$, converges for $|x|<1$ when $p=q+1$,
diverges for $p>q+1$ when $x\neq 0$ and the members of the series do not vanish.

Below we suggest that the differences  $a_i-a_j$,  $a_i-b_j$ are not integers.
\subsection{An ordinary differential equation}
\label{odupa}
 Let  $p\leq q$, put  $$\sigma=q-p.$$  Consider an ordinary differential equation

\begin{equation}
\label{odu}
(-1)^{q-p}(\prod_{j=1}^p(x\frac{d}{dx}+\mu_j)-\prod_{j=1}^q(x\frac{d}{dx}+\nu_j-1))G(x)=0.
\end{equation}

This equation in the case $p<q$ has two singular points: a regular point
 $0$ and an irregular point $\infty$ (in the case $p=q$ the singularities are the same but they are both regular).
 A base in the solution space can be described as follows. For $m\leq q$, $n\leq p $ introduce a $G$-function

\begin{align}
\begin{split}
&G_{p,q}^{m,n}(x;a_1,...,a_p,b_1,...,b_q)=\frac{1}{2\pi i}\int _L
\frac{\prod_{j=1}^m\Gamma(b_j-s)\prod_{j=1}^n\Gamma(1-a_j+s)}{\prod_{j=m+1}^q\Gamma(1-b_j+s)\prod_{j=n+1}^p\Gamma(a_j-s)}x^sds,
\end{split}
\end{align}

where the contour $L$  can be chosen for example as  follows (there are other natural choices):  it begins and ends in  $-\infty$,
encompasses  all poles of functions  $\Gamma(1-a_k+s)$ in positive direction, does not encompass the poles of functions
$\Gamma(b_j-s)$.  The resulting function can be continued until a multivalued analytic  function in $\mathbb{C}\setminus\{0\}$.

For the functions $G_{p,q}^{1,p}$ one has the following expansion in variable $x$:

\begin{align}
\begin{split}
\label{fg}
&G_{p,q}^{1,p}(x;a_1,...,a_p,b_1,...,b_q)=\frac{\prod_{j=1}^p\Gamma(1+b_1-a_j)}{\prod_{j=2}^q\Gamma(1+b_1-b_j)}\cdot\\&
\cdot x^{b_1}
F_{p,q-1}(1+b_1-a_1,...,1+b_1-a_p,1+b_1-b_2,...,1+b_1-b_q,-x).
\end{split}
\end{align}

The base of solutions of  \eqref{odu} in a neighborhood of  $0$ is
formed by functions

\begin{equation}
G_{p,q}^{1,p}(x;a_1,...,a_p,b_h,b_1,...,\widehat{b_h},...,b_q),\,\,\, h=1,...,q.
\end{equation}

where the parameters $\mu,\nu$ and  $a,b$ are related as follows

\begin{align}
\begin{split}
&a_j=1-\mu_j,\,\,\, j=1,...,p\\
&b_j=1-\nu_j,\,\,\, j=1,...,q.
\end{split}
\end{align}

There exists another base in the space of solutions.




\begin{align}
\begin{split}
\label{bif}
&G_{p,q}^{q,1}(x e^{\pi i};a_h,a_1,...,\widehat{a_h},...,a_p,b_1,...,b_{q}),\,\,\, h=1,..,p\\
& e^{\frac{2ih\lambda\pi }{\sigma}}G_{p,q}^{q,0}(xe^{2h\pi i};a_1,...,a_p,b_1,...,b_{q}),
\end{split}
\end{align}
where $\lambda$ is defined in \eqref{lda} and in the second formula $h\in I_{\sigma}$ and

\begin{equation}
I_{\sigma}=\{1-\frac{\sigma}{2},...,\frac{\sigma}{2}\}\text{ if $\sigma$ is even, } I_{\sigma}=\{-\frac{\sigma-1}{2},...,\frac{\sigma-1}{2}\}\text{ if $\sigma$ is odd }.
\end{equation}

Their asymptotic expansion is derived from the following result.
In the sector

\begin{equation}
|argx|<\pi(\frac{\sigma}{2}+1)
\end{equation}


  one has an asymptotic expansion for $x\rightarrow \infty$
\begin{align}
\begin{split}
\label{r1az1}
&G_{p,q}^{q,1}(x;a_h,a_1,...,\widehat{a_h},...,a_p,b_1,...,b_{q})\sim
\frac{x^{a_h-1}\Gamma(1+b_{q}-a_h)}{\Gamma(1+a_p-a_h)}\cdot\\&\cdot
F_{q,p-1}(1+b_1-a_h,...,1+b_q-a_h;1+a_1-a_h,...,1+a_p-a_h,-\frac{1}{x})
\end{split}
\end{align}

In the sector

\begin{equation}
|argx|<\pi(\sigma+\varepsilon),
\end{equation}

where $\varepsilon=\frac{1}{2}$ in the case $\sigma=1$ and $\varepsilon=1$ in the case $\sigma>1$,
one has an expansion


\begin{align}
\begin{split}
\label{r0az1}
&G_{p,q}^{q,0}(x;a_1,...,a_p,b_1,...,b_{q})\sim\frac{(2\pi)^{\frac{\sigma-1}{2}}}{\sigma^{\frac{1}{2}}}e^{-\sigma
x^{\frac{1}{\sigma}}}x^{\frac{\lambda}{\sigma}}M(\frac{1}{x^{\frac{1}{\sigma}}}),
\end{split}
\end{align}

where $M$ is a Taylor expansion in one variable  and a number $\lambda$
is defined in \eqref{lda} in the next Section.

Thus the set of exponential factors for the equation
\eqref{odu} is the following

\begin{equation}
\Phi=\{0,-\sigma (e^{2h\pi i}x)^{\frac{1}{\sigma}}\text{, where }
h\in I_{\sigma} \}.
\end{equation}

\subsection{The Stokes phenomenon and the monodromy for the equation \eqref{odu}}

The Stokes phenomenon and the monodromy for the equation  \eqref{odu}
are described in \cite{M}. Let us give this description.

\subsubsection{The formal monodromy at $\infty$}
\label{monod} A matrix $M$
of formal monodromy at  $\infty$ in the base \eqref{bif} is the following

\begin{enumerate}
\item If $\sigma=1$, then $M=diag(e^{-2i\pi
\mu_1},...,e^{-2i\pi\mu_p},e^{2 i \pi\lambda})$.
\item If $\sigma \geq 2$, then $M$ is block-diagonal, $M=diag(A,R)$, where
$$A=diag(e^{-2i\pi \mu_1},...,e^{-2i\pi\mu_p}),\,\, R=\begin{pmatrix}0
& 0&... &0& 1 \\  1 & 0 & ...& 0& 0\\ 0 & 1 & ... & 0 & 0\\
 & &...& &  \\ 0 & 0 & ...  & 1 & 0 \end{pmatrix}$$

\end{enumerate}

\subsubsection{The Stokes phenomenon}

From the formulas \eqref{r1az1},  \eqref{r0az1}  one obtains that $\infty$
is a ramified irregular singularity for \eqref{odu}. To remove the ramification one must change the independent variable

\begin{equation}
\label{raz2} t=x^{\frac{1}{\sigma}}.
\end{equation}

Following \cite{Du}  we give a description of the Stokes phenomenon in variable $x$.

The Stokes lines of  \eqref{odu} are the following. One series of Stokes lines is obtained if one takes $0$ and one of the exponential factors $-\sigma (e^{2h\pi
i}x)^{\frac{1}{\sigma}}$, these Stokes lines are defined by formula

\begin{equation}
\label{sl1} \{argt=-\frac{2h\pi}{\sigma}-\frac{\pi}{2}+n\pi,\,\,\, h\in I_{\sigma},\,\,\, n\in\mathbb{Z}\}.
\end{equation}

Also in the case  $\sigma>2$  one has a series of Stokes lines
corresponding to a pair of  exponential factors $-\sigma
(e^{2h\pi i}x)^{\frac{1}{\sigma}}$ for different $h$, these Stokes lines are defined by formula

\begin{equation}
\label{sl2} \{argt=\frac{(h_1+h_2)\pi}{\sigma}+n\pi,\,\,\, h_1,h_2\in I_{\sigma},\,\,\, h_1\neq h_2,\,\,\, n\in\mathbb{Z}\}.
\end{equation}

The Stokes sectors are defined by formulas

\begin{enumerate}
\item If $\sigma=1$, $2$, then the sectors are
\begin{equation}
\label{ss1}
 \Theta_n=\{x:  -\frac{\pi}{2}+(n-1)\pi<argx<
\frac{\pi}{2}+\pi n\}, \,\,\, n=0,1,2.
\end{equation}
\item If $\sigma\geq 3$ is odd then the sectors are
\begin{equation}
\label{ss2}
 \Theta_n=\{x:
-\frac{\pi}{2}+\frac{(n-1)\pi}{2\sigma}<arg x<
\frac{\pi}{2}+\frac{\pi n}{2\sigma}\},\,\,\,n=0,...,4\sigma
\end{equation}
\item  If $\sigma\geq 4$ is even then the sectors are
\begin{equation}
\label{ss3}
 \Theta_n=\{x:-\frac{\pi}{2}+\frac{(n-1)\pi}{2\sigma}<arg
x< \frac{\pi}{2}+\frac{\pi n}{2\sigma}\}, \,\,\, n=0,...,2\sigma.
\end{equation}
\end{enumerate}

The Stokes matrix corresponding to an intersection of $\Theta_n$
 and $\Theta_{n+1}$ is denoted as $S_n$. Let $\Sigma_n$,
 $\Sigma_{n+1}$ be solutions of \eqref{odu} that have an expansion \eqref{r1az1},  \eqref{r0az1}. Then for $\Theta_n\cap \Theta_{n+1}$

 \begin{equation}
\Sigma_{n}=\Sigma_{n+1}S_n.
 \end{equation}

Since a ramification takes place  relations between Stokes matrices $S_n$ take place.

For $\sigma=2$ $$S_1=M^{-1} S_0 M,$$

for $\sigma\geq 3$ let $m$ and $r$ be a quotient and a residue of division of $n$ by  $4$ in the case of odd $\sigma$ and of division by $2$ in the case of even  $\sigma$ $$S_n=M^{-m}S_r M^m.$$

 Instead of $S_n$ we sometimes write

\begin{enumerate}
\item  $S_{n\pi \setminus 2\sigma}$ in the case of odd $\sigma\geq
3$,
\item $S_{n\pi \setminus \sigma}$ in the case of even $\sigma\geq
4$,

 \item $S_{n\pi} $ in the case
$\sigma=1$, $2$.
\end{enumerate}

Thus we have the following independent Stokes matrices
\begin{enumerate}
\item in the case $\sigma=1$: $S_0$, $S_{\pi}$.
\item in the case $\sigma=2$: $S_0$.
\item in the case $\sigma\geq 3$: $S_0$, $S_{\pi \setminus
2\sigma}$, $S_{\pi \setminus \sigma}$, $S_{3\pi \setminus 2\sigma}$.
\item in the case $\sigma\geq 4$: $S_0$, $S_{\pi \setminus
\sigma}$.
\end{enumerate}

Below we use notation

 \begin{align}
 \begin{split}
 \label{lda}
 &\lambda=\frac{1}{2}(\sigma+1)+\sum_{i=1}^p\mu_i-\sum_{i=1}^q\nu_i,\\
&A_h:=-\frac{d^h}{dx^h}(\frac{\prod_{i=1}^{q}(1-xe^{2i\pi \nu_i})
}{\prod_{i=1}^{p}(1-xe^{2i\pi \mu_i})})\mid_{x=0},\\
 &B_h:=e^{2\pi i
\lambda}\frac{d^h}{dx^h}(\frac{\prod_{i=1}^{q}(1-xe^{2i\pi \nu_i})
}{\prod_{i=1}^{p}(1-xe^{2i\pi \mu_i})})\mid_{x=0}
,\\
&\zeta=e^{\frac{2\pi i}{\sigma}}..
\end{split}
 \end{align}

 \subsubsection{The Stokes matrices for $\sigma=1$}

The Stokes matrices are the following

 \begin{align*}
& S_0=I+\sum_{k=1}^{q-1}
2i\pi\frac{\prod_{i=1}^p\Gamma(1+\mu_k-\mu_i))}{\prod_{i=1}^q\Gamma(1+\mu_k-\nu_i))}E_{q,k},\\
& S_{\pi}=I+\sum_{k=1}^{q-1} 2i\pi e^{i\pi(\lambda+\mu_k)}
\frac{\prod_{i=1,i\neq
k}^p\Gamma(\mu_i-\mu_k)}{\Pi_{i=1}^q\Gamma(\nu_i-\mu_k))}E_{k,q}
 \end{align*}

 \subsubsection{The Stokes matrices for $\sigma=2$}

The Stokes matrices are the following

\begin{align*}
& S_0=I+A_1E_{q-1,q}+\sum_{k=1}^{q-2}4i\pi^2\frac{\prod_{i=1,i\neq
k}^p\Gamma(\mu_i-\mu_k)}{\prod_{i=1}^q\Gamma(\nu_i-\mu_k)}E_{k,q}+\\
&+\sum_{k=1}^{q-2}2i\pi\frac{\prod_{i=1}^p\Gamma(1+\mu_k-\mu_i)}{\prod_{i=1}^q\Gamma(1+\mu_k-\nu_i)}E_{q,k}
 \end{align*}

\subsubsection{The case $\sigma\equiv 0$ mod $4$, $\sigma\geq
4$}

\begin{align*}
& S_0=I+\sum_{k=1}^p(2\pi)^{\sigma}i\frac{\prod_{i=1,i\neq
k}^p\Gamma(\mu_i-\mu_k)}{\prod_{i=1}^q\Gamma(\nu_i-\mu_k)}E_{k,q}+\sum_{k=1}^p2\pi
\frac{\prod_{i=1}^p\Gamma(1+\mu_k-\mu_i)}{\prod_{i=1}^q\Gamma(1+\mu_k-\nu_i)}
E_{\sigma/2+p,k}+\\
&+\sum_{k=0}^{\sigma/4}\zeta^{-2\lambda k} A_{2k}
E_{3\sigma/4+p-k,3\sigma/4+p+k}+\sum_{k=1}^{\sigma/4-1}e^{-2i\pi\lambda}\zeta^{2\lambda
k} B_{2k}E_{\sigma/4+p+k,\sigma/4+p-k},\\
&S_{\pi/\sigma}=I+\sum_{k=1}^{\sigma/4}\zeta^{-\lambda(2k-1)}A_{2k-1}E_{3\sigma/4+p-k,3\sigma/4+p+k-1}+\\&+\sum_{k=1}^{\sigma/4-1}
e^{-2i\pi\lambda}\zeta^{\lambda(2k-1)}B_{2k-1}E_{\sigma/4+p+k-1,\sigma/4+p-k}+e^{-i\lambda\pi}
B_{\sigma/2-1}E_{\sigma/2+p-1,q}.
 \end{align*}

\subsubsection{The case $\sigma\equiv 2$ mod $4$}

\begin{align*}
& S_0=I+\sum_{k=1}^p(2\pi)^{\sigma}i
2\pi)^{\sigma}i\frac{\prod_{i=1,i\neq
k}^p\Gamma(\mu_i-\mu_k)}{\prod_{i=1}^q\Gamma(\nu_i-\mu_k)}E_{k,q}+\sum_{k=1}^p2\pi
\frac{\prod_{i=1}^p\Gamma(1+\mu_k-\mu_i)}{\Pi_{i=1}^q\Gamma(1+\mu_k-\nu_i)}
E_{\sigma/2+p,k}+\\
&+\sum_{k=1}^{(\sigma+2)/4}\zeta^{-\lambda(2k-1)}A_{2k-1}E_{3(\sigma+2)/4+p-k,3(\sigma-2)/4+p+k+1}+\\
&+\sum_{k=1}^{(\sigma-2)/4}e^{-2i\pi\lambda}\zeta^{\lambda(2k-1)}B_{2k-1}E_{(\sigma-2)/4+p+k-1,(\sigma+2)/4+p-k},\\
&S_{\pi/\sigma}=I+\sum_{k=1}^{(\sigma-2)/4}\zeta^{-2\lambda
k}A_{2k}E_{3(\sigma-2)/4+p-k,3(\sigma-2)/4+p+k+1}+\\
&+e^{-i\pi\lambda}\zeta^{-\lambda}
B_{\sigma/2-1}E_{\sigma/2+p-1,q}+\sum_{k=1}^{(\sigma-2)/4-1}e^{-2i\pi\lambda}\zeta^{2\lambda
k}B_{2k}E_{(\sigma-2)/4+p+k-1,(\sigma-2)/4+p-k}.
\end{align*}

\subsubsection{The case of odd $\sigma$ }

\begin{align*}
& S_0=I+\sum_{k=1}^p2i\pi
\frac{\prod_{i=1}^p\Gamma(1+\mu_k-\mu_i))}{\prod_{i=1}^q\Gamma(1+\mu_k-\nu_i))}
E_{(p+q+1)/2,k},\\
&S_{\pi/\sigma}=I+\sum_{k=1}^p (2\pi )^{\sigma} e^{i\pi
(\mu_k+(\lambda/\sigma))} \frac{\prod_{i=1,i\neq
k}^p\Gamma(\mu_i-\mu_k)}{\prod_{i=1}^q\Gamma(\nu_i-\mu_k)}E_{k,q}.
\end{align*}

\subsubsection{The case $\sigma\equiv 1$ mod $4$}

\begin{align*}
&S_{\pi/ 2\sigma}=I+\sum_{k=1}^{(\sigma-1)/4}\zeta ^{-2\lambda k}
A_{2k}E_{3(\sigma-1)/4+p-k+1,3(\sigma-1)/4+p+k+1}+\\
&+\sum_{k=1}^{(\sigma-1)/4} e^{-2\pi i \lambda}
\zeta^{\lambda(2k-1)}B_{2k-1}
E_{(\sigma-1)/4+p+k,(\sigma-1)/4+p-k+1},\\&
S_{3\pi/2\sigma}=I+\sum_{k=1}^{(\sigma-1)/4}\zeta^{-\lambda(2k-1)}A_{2k-1}E_{3(\sigma-1)/4+p-k+1,3(\sigma-1)/4+p+k}+\\
&+\sum_{k=1}^{(\sigma-1)/4-1}e^{-2i\pi \lambda} \zeta^{2\lambda k}
B_{2k}E_{(\sigma-1)/4+p+k,(\sigma-1)/4+p-k}+ e^{-i\lambda
\pi}\zeta^{-\lambda/2} B_{(\sigma-1)/2} E_{(\sigma-1)/2+p,q}.
\end{align*}

\subsubsection{The case $\sigma\equiv 3$ mod $4$}

\begin{align*}
& S_{\pi/
2\sigma}=I+\sum_{k=1}^{(\sigma+1)/4}\zeta^{-\lambda(2k-1)}A_{2k-1}E_{3(\sigma+1)/4+p-k,3(\sigma+1)/4+p+k-1}+\\
&+\sum_{k=1}^{(\sigma+1)/4-1} e^{-2\pi i \lambda} \zeta^{2\lambda
k}B_{2k}
E_{(\sigma+1)/4+p+k,(\sigma+1)/4+p-k},\\
&S_{3\pi/2\sigma}=I+\sum_{k=1}^{(\sigma+1)/4-1}\zeta ^{-2\lambda k}
A_{2k}E_{3(\sigma+1)/4+p-k-1,3(\sigma+1)/4+p+k-1}+\\
&+\sum_{k=1}^{(\sigma+1)/4-1}e^{-2i\pi \lambda} \zeta^{\lambda
(2k-1)} B_{2k-1}E_{(\sigma+1)/4+p+k-1,(\sigma+1)/4+p-k}+
e^{-i\lambda \pi}\zeta^{-\lambda/2} B_{(\sigma-1)/2}
E_{(\sigma+1)/2+p,q}.
\end{align*}


\section{The Gelfand-Kapranov-Zelevinsky system}

\label{sgkz}

The Gelfand-Kapranov-Zelevinsky system (the GKZ system, see \cite{Gel}, \cite{GG})
is defined as follows. Let $B\subset \mathbb{Z}^n$ be a lattice.
Then a GKZ system associated with $B$  is a system of partial differential equations

\begin{align}
\begin{split}
\label{gkz0}
&\sum_{i=1}^Na_i x_i\frac{\partial}{\partial x_i}y-\beta y=0,\,\,\, a=(a_1,...,a_n)\in annL\\
&(\frac{\partial}{\partial x_i})^cy=(\frac{\partial}{\partial
x_i})^{c'}y\,\,\,\Leftrightarrow c-c'\,\,\,\in B,
\end{split}
\end{align}
where $L$ is a subspace spanned by $B$, and  $\beta$ is presented as $<\gamma,a>$ for some $\gamma\in \mathbb{C}^n$.

To a system of PDE there corresponds a  $\mathcal{D}$-module (see \cite{dmo}). In our case this module is holonomic (see \cite{A}) and regular if and only is $(1,...,1)\perp B$. The holonomicity implies that this $\mathcal{D}$-module on a dense open subset is a  $\mathcal{D}$-module of sections of a bundle with a connection.  Thus one can apply to this system the theory above.

In the present paper we consider a special case of a GKZ system. Let $n=p+q+1$ and let $e_i=(0,...,1_{i},...,0)$ denote a standard base vector in $\mathbb{C}^n$. Define the lattice $B$
\begin{equation*}
B=\mathbb{Z}(-(e_1+...+e_p)+(e_{p+1}+...+e_{p+q+1})),
\end{equation*}

Without less of generality we suggest that $p\leq q+1$.
Choose a base in $annL$

\begin{align*}
&a_{12}=(1,-1,...,0_{p},0_{p+1},...,0),\\
&...\\
&a_{1,p}=(1,0,...,-1_{p},0_{p+1},...,0),\\
&a_{1,p+1}=(1,0,...,0_{p},1_{p+1},...,0),\\
&...\\
&a_{1,p+q+1}=(1,0,...,0_{p},0_{p+1},...,1).\\
\end{align*}

Then \eqref{gkz0} turns into
\begin{align}
\begin{split}
\label{gkzb}
&x_1\frac{\partial y}{\partial x_1}-x_2\frac{\partial y}{\partial x_2}-\beta_{1,2} y=0,\\
&...\\
&x_1\frac{\partial y}{\partial x_1}-x_{p}\frac{\partial y}{\partial x_{p}}-\beta_{1,p} y=0,\\
&x_1\frac{\partial y}{\partial x_1}+x_{p+1}\frac{\partial y}{\partial x_{p+1}}-\beta_{1,p+1} y=0,\\
&...\\
&x_1\frac{\partial y}{\partial x_1}+x_{p+q+1}\frac{\partial y}{\partial x_{p+q+1}}-\beta_{1,p+q+1} y=0,\\
&\frac{\partial^p y}{\partial x_1...\partial x_p}=\frac{\partial^{q+1}
y}{\partial x_{p+1}...\partial x_{p+q+1}},
\end{split}
\end{align}

Also let  $\gamma=(-a_1,...,-a_p,b_1-1,...,b_q-1,0)$. Then the parameters  $\beta_{1,i} $ have the following presentation $\beta_{1,2}=((-a_1)-(-a_2))$,..., $\beta_{1p}=((-a_1)-(-a_p))$,
$\beta_{1,p+1}=((-a_1)+(b_1-1))$,...,$\beta_{1,p+q+1}=((-a_1)+(b_q-1))
$.

The system is regular if and only if $p=q+1$. We suggest below that $$p<q+1.$$

\section{The singular divisor}

The characteristic variety of the system  \eqref{gkzb} is a submanifold in $T^*X$,defined by equations
\begin{align}
\begin{split}
&x_1\xi_1=x_2\xi_2,\,\,...,\,\,x_1\xi_1=x_p\xi_p,\\
&x_1\xi_1=-x_{p+1}\xi_{p+1},\,\,...,\,\,x_1\xi_1=x_{p+q+1}\xi_{p+q+1},\\
&\xi_1...\xi_p=0\text{ if $p>q+1$},\xi_{p+1}...\xi_{p+q+1}=0\text{
if $q+1>p$},\\&\xi_1...\xi_p=\xi_{p+1}...\xi_{p+q+1}\text{ if
$p=q+1$}
\end{split}
\end{align}

Remove from this submanifold a zero section onto coordinates $x_1,...,x_{p+q+1}$.

In the case  $p<q+1$ one obtains a manifold  defined be equations

\begin{equation}
x_{1}...x_{p}=0.
\end{equation}

In the case $p=q+1$, that is in the regular case  one obtains a manifold  defined be equations

\begin{equation}
x_1...x_p-x_{p+1}...x_{p+q+1}=0.
\end{equation}

\section{ The space $\widetilde{X}$ for the system \eqref{gkzb}}
Following the Section \ref{pvo} let us give a description of the space   $\widetilde{X}$ for the considered GKZ system.
 We consider only the case $p<q+1$.  Then the singular divisor $D$ is
 $$D=D_1+...+D_p,\,\,\,\, D_i=\{x_i=0\}.$$

Introduce on $\mathbb{C}^{p+q+1}$ bundles $L(D_1)\simeq
L\{0\}\widehat{\otimes} \mathbb{C}\widehat{\otimes}
...\widehat{\otimes} \mathbb{C}, $ where $L\{0\}$ is a linear bundle $\mathbb{C}$ associated with the divisor $\{0\}$. Also construct analogous bundles
$L(D_2)$,...,$L(D_p)$.

For the stalks one has an isomorphism
\begin{equation}
S^1(L\{0\}\widehat{\otimes} \mathbb{C}\widehat{\otimes}
...\widehat{\otimes} \mathbb{C})_{(x_1,...,x_{p+q+1})}\simeq S^1
L\{0\}_{x_1}.
\end{equation}

Thus, one has

\begin{equation}
L:=S^1(L(D_1)\oplus...\oplus L(D_p))_{(x_1,...,x_{p+q+1})}=
S^1L\{0\}_{x_1}\oplus...\oplus S^1 L\{0\}_{x_p}.
\end{equation}

On the set  $X^*=X\setminus D$ there is a mapping $X^*\rightarrow L$ that sends $x$ to
$(\frac{x_1}{|x_1|},...,\frac{x_p}{|x_p|})=(arg x_1,...,arg x_{p}).$
The closure of  $L$ of an image of this mapping is $\widetilde{X}$.

Thus a subset of  $\widetilde{X}$ that is mapped to
$D_{i_1}\cap...\cap D_{i_k}\setminus(\bigcup_{i\neq
i_j}D_i)=\{x_{i_1}=...=x_{i_k}=0,
x_1,...,\widehat{x_{i_1}},...\widehat{x_{i_k}},...,x_p\neq 0\}$ is
homeomorphic to $\underbrace{S^{1}\times...\times S^{1}}_{p \text{
times }}.$

\section{Solutions of the system \eqref{gkzb}}
To describe the Stokes phenomenon we need an explicit description of the space of solutions of the GKZ system, let us give it.
\subsection{$\Gamma$-series}

The right parts  $\beta_{1,i}$ in the system \eqref{gkzb} can be represented as
$(a_{1,i},\gamma)$, where $\gamma=(-a_1,...,-a_p,b_1-1,...,b_q-1,0)$.

There exists a base in the space of solution of the system \eqref{gkzb} whose elements are written as $\Gamma$-series.

\begin{defn}
{\it A $\Gamma$-series} associated with a lattice $B$ and a vector $\gamma$
is defined by the formula

\begin{equation}
F_{B}(\gamma,x)=\sum_{b\in
B}\frac{x^{\gamma+b}}{\Gamma(\gamma+b+1)},
\end{equation}
where $x^{\gamma+b}:=x_1^{\gamma_1+b_1}...x_N^{\gamma_N+b_N}$ и
$\Gamma(\gamma+b+1):=\Gamma(\gamma_1+b_1+1)...\Gamma(\gamma_N+b_N+1)$
\end{defn}

A space of formal solutions of the system \eqref{gkzb}is generated by series $F_{B}(\gamma^1,x),...,F_{B}(\gamma^q,x)$, where $\gamma^i_{p+i}\in \mathbb{Z}$ and
$(\gamma^i,a_{1j})=\beta_{1j}$, $j=1,...,p+q+1$. The only relations between these series are
$$F_{B}(\gamma+b,x)=F_{B}(\gamma,x),\,\,\, b\in B.$$

The equations $(\gamma^i, a_{1,j})=\beta_{1,j}$ are written explicitly as follows

\begin{align*}
&\gamma_1-\gamma_2=(-a_1)-(-a_2),\\
&...\\
&\gamma_1-\gamma_p=(-a_1)-(-a_p),\\
&\gamma_1+\gamma_{p+1}=(-a_1)+(b_1-1),\\
&...\\
&\gamma_1+\gamma_{p+q+1}=(-a_1)+(b_{q+1}-1).
\end{align*}

Put $\gamma_{p+i}=0$, then one has

\begin{align*}
&\gamma_1=(-a_1)+(b_i-1),\\
&\gamma_2=(-a_2)+(b_i-1),\\
&...\\
&\gamma_p=(-a_p)+(b_i-1),\\
&\gamma_{p+1}=-(b_i-1)+(b_1-1),\\
&...\\
&\gamma_{p+q+1}=-(b_i-1)+(b_{q+1}-1).
\end{align*}

Take as base solutions the series for which $\gamma^i_{p+i}=0$ explicitly the vectors  $\gamma$ are written as follows

\begin{align}
\begin{split}
&\gamma_i=((-a_1)+(b_{i}-1),...,(-a_p)+(b_{i}-1),b_{1}-b_{i},...,0_{p+i},...,b_{q+1}-b_{i}),\\
&i=1,...,q+1.
\end{split}
\end{align}

The considered  $\Gamma$-series ar related with the hypergeometric series of one argument $F_{p,q}$ by the equality

\begin{align*}
& F_{B}(-a_1,...,-a_p,b_1-1,...,b_q-1,0;x_1,...,x_{p+q+1})=\\
&=
\frac{1}{\Gamma(1-a_1)...\Gamma(1-a_p)\Gamma(b_1)...\Gamma(b_q)}x_1^{-a_1}...x_p^{-a_p}
x_{p+1}^{b_1-1}...x_{p+q}^{b_q-1}\\&
F_{p,q}(a_1,...,a_p;b_1,...,b_q;\frac{x_{p+1}...x_{p+q+1}}{x_1...x_p}).
\end{align*}

From this relation one sees that   $\Gamma$-series converges for $x_1,...,x_p\neq 0$.  Thus we have a description not of the space of formal solutions
 but of the space of analytic solutions of  \eqref{gkzb}.

Using these relation the solutions of \eqref{gkzb} can be written as follows

\begin{align}
\begin{split}
\label{fpqsol}
&F_B(\gamma_i,x_1,...,x_{p+q+1})=\\&= F_{B}(b_{i}-a_1-1,...,b_{i}-a_p-1,b_{1}-b_{i},...,0_{i},...,b_{q+1}-b_{i};x_1,...,x_{p+q+1})=\\
&=\frac{1}{\Gamma(b_{i}-a_1)...\Gamma(b_{i}-a_p)\Gamma(b_{1}-b_{i}+1)
...\Gamma(b_{q+1}-b_{i}+1)}x_1^{b_{i}-a_1-1}...x_p^{b_{i}-a_p-1}x_{p+1}^{b_{1}-b_{i}}
...\widehat{x_i}...x_{p+q+1}^{b_{q+1}-b_{i}}\\&
F_{p,q}(a_{1}-b_i+1,...,a_{p}-b_i+1,b_{1}-b_{i}+1,...,b_{q+1}-b_{i}+1;\frac{x_{p+1}...x_{p+q+1}}{x_1...x_p}).
\end{split}
\end{align}

\section{An isomorphism between the spaces of solution of the equation \eqref{odu} and of the system \eqref{gkzb}}

Using \eqref{fg}, the solutions \eqref{fpqsol} can be written as follows

\begin{align}
\begin{split}
&\frac{\prod_{j=1,j\neq i}^{q+1
}\Gamma(1+(-b_i)-(-b_j))}{\prod_{j=1}^p\Gamma(1+(-b_i)-(-a_j))\prod_{j=1}^p\Gamma(b_i-a_j)\prod_{j=1,j\neq
i}^{q+1}
\Gamma((-b_i)+b_j+1)}\cdot\\
&\cdot \prod_{j=1}^px_j^{b_i-a_j-1}\prod_{j=1}^{q+1}x_{j+p}^{b_j-b_i}(\frac{x_{p+1}...x_{p+q+1}}{x_1...x_p})^{b_i}\cdot\\
&\cdot G_{p, q+1}^{1,p}(-a_1,...,-a_p,-b_i,-b_1,...,\widehat{-b_i},...,-b_{q+1},-\frac{x_{p+1}...x_{p+q+1}}{x_1...x_p})=\\
&=\frac{1}{\prod_{j=1}^p\Gamma(1+(-b_i)-(-a_j))\prod_{j=1}^p\Gamma(b_i-a_j)} \prod_{j=1}^px_j^{-a_j-1}\prod_{j=1}^{q+1}x_{j+p}^{b_j}\cdot\\
&\cdot G_{p,
q+1}^{1,p}(-a_1,...,-a_p,-b_i,-b_1,...,\widehat{-b_i},...,-b_{q+1},-\frac{x_{p+1}...x_{p+q+1}}{x_1...x_p})
\end{split}
\end{align}

Thus we have proved the following

\begin{thm}

There exists an isomorphism
\begin{equation}
\label{izom} f(x)\mapsto
\frac{1}{\prod_{j=1}^p\Gamma(1+(-b_i)-(-a_j))\prod_{j=1}^p\Gamma(b_i-a_j)}
\prod_{j=1}^px_j^{-a_j-1}\prod_{j=1}^{q+1}x_{j+p}^{b_j}
f(-\frac{x_{p+1}...x_{p+q+1}}{x_1...x_p}),
\end{equation}

the maps the  analytic solutions of the equation \eqref{odu} to the solutions of the system \eqref{gkzb}.

\end{thm}

\section{The Riemann-Hilbert functor}
\label{stoxgkz}

To describe the Stokes phenomenon it is more useful to use not the
base $F_B(\gamma_i,x)$, but a base corresponding to \eqref{bif}
under the isomorphism \eqref{izom}. It consists of functions with
the following asymptotics.

Define a sector and integers $k,g$ using formulas
\eqref{secto} in which we substitute
$x=\frac{x_{p+1}...x_{p+q+1}}{x_1...x_p}$.

Then using the correspondence \eqref{izom} we get that there exist a base in the space of solutions of the GKZ system that has the following asymptotic
 in the sector described above.

At first there are function for $h=1,...,p$ that have the asymptotic

\begin{align}
\begin{split}
\label{raze1}
&\frac{1}{\prod_{j=1}^p\Gamma(1+(-b_i)-(-a_j))\prod_{j=1}^p\Gamma(b_i-a_j)}
\prod_{i=1}^px_i^{-a_1-1}\prod_{i=1}^{q+1}x_{i+p}^{b_1}\\&
(e^{i\pi
}\frac{x_{p+1}...x_{p+q+1}}{x_1...x_p})^{a_h-1}\frac{\Gamma(1+b_{q}-a_h)}{\Gamma(1+a_p-a_h)}\cdot\\&\cdot
F_{q+1,p-1}(1+b_1-a_h,...,1+b_{q+1}-a_h;1+a_1-a_h,...,1+a_p-a_h,-e^{i\pi}\frac{x_1...x_p}{x_{p+1}...x_{p+q+1}})
\end{split}
\end{align}

Secondly the are functions for $h\in I_{\sigma}$ that have asymptotic

\begin{align}
\begin{split}
\label{raze2}
&\frac{1}{\prod_{j=1}^p\Gamma(1+(-b_i)-(-a_j))\prod_{j=1}^p\Gamma(b_i-a_j)}
\prod_{i=1}^px_i^{-a_1-1}\prod_{i=1}^{q+1}x_{i+p}^{b_1}
\frac{(2\pi)^{\frac{\sigma-1}{2}}}{\sigma^{\frac{1}{2}}}\\&exp(-\sigma
(-e^{2h\pi
i}\frac{x_{p+1}...x_{p+q+1}}{x_1...x_p})^{\frac{1}{\sigma}})(-\frac{x_{p+1}...x_{p+q+1}}{x_1...x_p})^{\frac{\lambda}{\sigma}}M(\frac{1}{(-e^{2h\pi
i}\frac{x_{p+1}...x_{p+q+1}}{x_1...x_p})^{\frac{1}{\sigma}}}),
\end{split}
\end{align}

 where  $M$ is a Taylor series of the argument written in the formula.



\subsection{The set of exponential factors}
Although it is hard to write the Leveltt-Turittin decomposition for the GKZ system, using formulas  \eqref{raze1},
\eqref{raze2} one can easily write the set of exponential factors in this decomposition.



More precise the formulas   \eqref{raze1}, \eqref{raze2} give the
following set of exponential factors
\begin{equation}
\label{expf}
 \Phi=\{0, (-\sigma)(-e^{2h\pi
i}\frac{x_{p+1}...x_{p+q+1}}{x_1...x_p})^{\frac{1}{\sigma}}\text{
 where  }h\in I_{\sigma}\}.
\end{equation}

Using the isomorphism \eqref{izom} one obtains the following. If  $$\Phi=\{\psi_i(x)\}$$  is a  set of exponential factors of the equation \eqref{odu}, the set of exponential factors of the
 system \eqref{gkzb} is given by
$$\Phi=\{\psi_i(-\frac{x_{p+1}...x_{p+q+1}}{x_1...x_p})\}.$$

In particular the system \eqref{gkzb} has a good formal
decomposition along a singular divisor.

\subsection{The sheaf $\mathcal{L}$}

Let us give a description of the solution sheaf  $\mathcal{L}$. This a locally trivial vector bundle on   $X^*$ of rank  $q+1$. The monodromy operator corresponding to a bypass along  $x_i=0$  for $i=1,...,p$ equals to $M_{\infty}$, where $M_{\infty}$  is a monodromy of the equation \eqref{odu}, which is expressed through the formal monodormy
 $M$ described in  \ref{monod} and Stokes matrices as follows (see \cite{I}):
$$M_{\infty}=M...S_1S_0.$$

\subsection{The sheaf $\mathcal{L}_{\leq}$}

Define a sheaf $\mathcal{L}_{\leq}$ on the space $\mathcal{J}^{et}$.
Let $x\in D_{i_1}\cap...\cap D_{i_l}$, where
$i_1,...,i_l\in\{1,...,p\}$, но при этом $x\notin
D_1\cup...\widehat{D_{i_1}}\cup...\cup\widehat{ D_{i_l}}\cup...\cup
D_p$.  Without less of generality we suggest that  $x\in
(D_1\cap...\cap D_l)\setminus(D_{l+1}\cup...\cup D_p)$

\subsubsection{Stokes hypersurfaces and hypersectors}

Since the set of exponential factors is given by
\eqref{expf} we have Stokes hypersurfaces or type

\begin{equation}
\label{sps1} St(0,-\sigma(-e^{2h\pi
i}\frac{x_{p+1}...x_{p+q+1}}{x_1...x_p})^{\frac{1}{\sigma}}),
\end{equation}

and for $\sigma>1$ also Stokes hypersurfaces of type

\begin{equation}
\label{sps2} St(-\sigma(-e^{2h_1\pi
i}\frac{x_{p+1}...x_{p+q+1}}{x_1...x_p})^{\frac{1}{\sigma}},-\sigma(-e^{2h_2\pi
i}\frac{x_{p+1}...x_{p+q+1}}{x_1...x_p})^{\frac{1}{\sigma}}).
\end{equation}

 On an intersection $(D_1\cap...\cap
D_l)\setminus(D_{l+1}\cup...\cup D_p)$ for the Stokes hypersurfaces
of type \eqref{sps1} we get
$u(0,...,0,x_{p+1},...,x_{p+q+1})=-\sigma(-e^{2h\pi
i}(\frac{x_{p+1}...x_{p+q+1}}{x_{l+1}...x_p})^{\frac{1}{\sigma}}$,
and the hypersurface is defined by formula  (see the formula \eqref{sl1})
\begin{align*}
&\frac{\pi}{\sigma}-\frac{2h}{\sigma}-\frac{arg x_{l+1}}{\sigma}-...-\frac{arg
x_p}{\sigma}+\frac{arg x_{p+1}}{\sigma}+...+\frac{arg
x_{p+q+1}}{\sigma}+\frac{\theta_1}{\sigma}+...+\frac{\theta_l}{\sigma}=\\&=\pm\frac{\pi}{2\sigma}+\frac{2\pi n}{\sigma}
\,\,\,\, n\in \mathbb{Z}
\end{align*}

In the case  $\sigma>1$ for the hypersurfaces of type \eqref{sps2}
we have $u(0,...0,x_{l+1},...,x_{p+q+1})=-\sigma(-e^{2h_1\pi i}+e^{2h_2\pi
i})^{\frac{1}{\sigma}}(\frac{x_{p+1}...x_{p+q+1}}{x_{l+1}...x_p})^{\frac{1}{\sigma}}$,
and the hypersurface itself is defined by the equality (see the formula \eqref{sl2})
\begin{align*}
&\frac{\pi}{\sigma}-\frac{(h_1+h_2)\pi}{\sigma^2}-\frac{arg x_{l+1}}{\sigma}-...-\frac{arg
x_p}{\sigma}+\frac{arg x_{p+1}}{\sigma}+...+\frac{arg
x_{p+q+1}}{\sigma}+\frac{\theta_1}{\sigma}+...+\frac{\theta_l}{\sigma}=\\&=\pm\frac{\pi}{2\sigma}+\frac{2\pi n}{\sigma}
\,\,\,\, n\in \mathbb{Z}
\end{align*}

Thus for fixed
$x_{l+1},...,x_{p+q+1},\theta_2,...,\theta_l$ on the circle
\begin{equation}
\label{okr} S^1\times \{\theta_2\}\times...\times
\{\theta_l\}\subset
\omega^{-1}(0,...,0,x_{l+1},...,x_{p+q+1})\subset
\partial \widetilde{X}
\end{equation}
 with a coordinate $\theta_i$ we have Stokes directions obtained from Stokes directions of the equations \eqref{odu} by rotating by the angle
\begin{equation}
\label{ugol} \frac{\pi}{\sigma}-\frac{arg x_{l+1}}{\sigma}-...-\frac{arg
x_p}{\sigma}+\frac{arg x_{p+1}}{\sigma}+...+\frac{arg
x_{p+q+1}}{\sigma}+\frac{\theta_2}{\sigma}+...+\frac{\theta_l}{\sigma},
\end{equation}

when one changes coordinates
$x_{l+1},...,x_{p+q+1},\theta_2,...,\theta_l$ on the torus
$S^1\times...\times S^1=\omega^{-1}(0,...,0,t_{l+1},...,t_{p+q+1})$
 these directions cover the Stokes hypersurfaces.

An analogous construction can be applied to the Stokes sectors of the equation  \eqref{odu}. On the circle \eqref{okr} one has Stokes sectors that are obtained from the Stokes sectors of the equation by rotating by the angle
\eqref{ugol}. If one changes the coordinates
$x_{l+1},...,x_{p+q+1},\theta_2,...,\theta_l$ these Stokes sectors are rotating and they cover the set which we call {\it the Stokes hypersector}.

\subsubsection{The local system $gr\mathcal{L}_{\psi}$}

Introduce a local system у $gr_0\mathcal{L}$ on the space
$\mathcal{J}^{et}$ of dimension $p$, and the systems
$gr_{-\sigma(-e^{2h\pi
i}\frac{x_{p+1}...x_{p+q+1}}{x_1...x_p})^{\frac{1}{\sigma}}}\mathcal{L}$, $h\in I_{\sigma}$, 
of dimension $1$.

\begin{defn} A sheaf $gr_0\mathcal{L}$  is an inverse image under the action of $\mu$ on
$\partial\widetilde{X}$ of a sheaf formed by solutions of the system
 \eqref{gkzb}, that have power-like growth along $(D_1\cap...\cap
D_l)\setminus(D_{l+1}\cup...\cup D_p)$.
\end{defn}

\begin{defn}
A sheaf $gr_{-\sigma(-e^{2h\pi
i}\frac{x_{p+1}...x_{p+q+1}}{x_1...x_p})^{\frac{1}{\sigma}}}\mathcal{L}$
 is an inverse image under the action of $\mu$ on the space
$\partial\widetilde{X}$ of a sheaf formed by solutions of the system
 \eqref{gkzb}, that  grow as $e^{-\sigma(-e^{2h\pi
i}\frac{x_{p+1}...x_{p+q+1}}{x_1...x_p})^{\frac{1}{\sigma}}}$  near
$(D_1\cap...\cap D_l)\setminus(D_{l+1}\cup...\cup D_p)$.
\end{defn}

\subsubsection{The gluing of  $\mathcal{L}_{\leq}$ over an intersection of Stokes hypersectors }

In Section \ref{rhf} it is stated that if there exist a good formal decomposition if one restricts  $\mathcal{L}_{\leq}$ onto an inverse image in $\mathcal{J}^{et}$ of  a Stokes hypersector  under the action of $\mu^{-1}$ one gets a decomposition
 \begin{equation}
 \label{razlll}
\mathcal{L}_{\leq \varphi(x)}=\bigoplus_{\psi\in\Phi, \,\,\,
\psi(x)\leq \varphi(x)}gr_{\psi(x)}\mathcal{L}\mid_{\varphi(x)}.
 \end{equation}

In inverse images of different hypersectors these decompositions are different.






To define the sheaf $\mathcal{L}_{\leq}$ one has to explain how to glue the decompositions \eqref{razlll} in an intersection of two inverse images of two Stokes hypersectors.
Due to the existence of an isomorphism \eqref{izom} we have the following.
An intersection of two Stokes hypersectors for the system \eqref{gkzb}
corresponds to an intersection of two Stokes sectors for the equation \eqref{odu}. Then  two decompositions \eqref{razlll}  in an intersection of two inverse images of two Stokes hypersectors
 are related by a corresponding Stokes matrix of the equation \eqref{odu}.

\end{document}